\title{\vspace{-0.6cm} Quasi-randomness of graph balanced cut properties}
\author{Hao Huang \thanks{Department of Mathematics, UCLA, Los
Angeles, CA, 90095. Email: {\tt huanghao@math.ucla.edu}.} \and Choongbum
Lee\thanks{Department of Mathematics, UCLA,  Los Angeles, CA 90095.
Email: {\tt abdesire@math.ucla.edu}.
Research supported in part by Samsung Scholarship. }}
\date{}

\documentclass[11pt]{article}
\usepackage{amsfonts,amssymb,amsmath,latexsym}
\usepackage{graphicx}
\newcommand{\PCUT}[1]{\mathcal{P}_{C}( #1 )}
\newcommand{\VALPHA}{\vec{\alpha}}   

\numberwithin{figure}{section}
\newtheorem{thm}{Theorem}[section]

\newtheorem{lemma}[thm]{Lemma}

\newtheorem{claim}[thm]{Claim}

\newtheorem{ques}[thm]{Question}

\newtheorem{dfn}[thm]{Definition}
\newenvironment{pf}
      {\medskip\noindent{\bf Proof.}\hspace{1mm}}
      {\hfill$\Box$\medskip}

\def\qed{\ifvmode\mbox{ }\else\unskip\fi\hskip 1em plus 10fill$\Box$}

\begin{document}
\maketitle

\abstract{Quasi-random graphs can be informally described as graphs
whose edge distribution closely resembles that of a truly random
graph of the same edge density. Recently, Shapira and
Yuster proved the following result on quasi-randomness
of graphs. Let $k \ge 2$ be a fixed integer, $\alpha_1, \ldots,
\alpha_k$ be positive reals satisfying $\sum_{i} \alpha_i = 1$ and
$(\alpha_1, \ldots, \alpha_k) \neq (1/k, \ldots, 1/k)$, and $G$ be a
graph on $n$ vertices. If for every partition of the vertices of $G$
into sets $V_1, \ldots, V_k$ of size $\alpha_1 n, \ldots, \alpha_k
n$, the number of complete graphs on $k$ vertices which have exactly
one vertex in each of these sets is similar to what we would expect
in a random graph, then the graph is quasi-random.
However, the method of quasi-random hypergraphs they used did
not provide enough information to resolve the case $(1/k, \ldots, 1/k)$ for graphs. In their work,
Shapira and Yuster asked whether this case also
forces the graph to be quasi-random. Janson also posed the same question
in his study of quasi-randomness under the framework of graph limits.
In this paper, we positively answer their question.
}

\section{Introduction} \label{introduction}

The study of random structures has seen a tremendous success in
modern combinatorics and theoretical computer science. One example
is the Erd\H{o}s-R\'enyi random graph $G(n,p)$ proposed in the
1950's and intensively studied thereafter. $G(n,p)$ is the
probability space of graphs over $n$ vertices where each pair of
vertices forms an edge independently with probability $p$. Random
graphs are not only an interesting object of study on their own but
also proved to be a powerful tool in solving numerous open problems.
The success of random structures served as a natural motivation for
the following question: How can one tell when a given structure
behaves like a random one? Such structures are called
\emph{quasi-random}. In this paper we study quasi-random graphs,
which, following Thomason \cite{Thomason1, Thomason2}, can be informally defined
as graphs whose edge distribution closely resembles that of a random graph (the
formal definition will be given later).
One fundamental result in the study of quasi-random graphs is the following
theorem proved by Chung, Graham and Wilson \cite{chung1989quasi} (here we only
state part of their result).



\begin{thm} \label{thm_chunggraham}
Fix a real $p \in (0,1)$. For an $n$-vertex graph $G$, define $e(U)$
to be the number of edges in the induced subgraph spanned by vertex
set $U$, then the following properties are equivalent.

$\mathcal{P}_1$: For any subset of vertices $U \subset V(G)$, we
have $e(U)=\frac{1}{2}p|U|^2 \pm o(n^2)$.

$\mathcal{P}_2(\alpha)$: For any subset of vertices $U \subset V(G)$
of size $\alpha n$, we have  $e(U)=\frac{1}{2}p|U|^2 \pm o(n^2)$.

$\mathcal{P}_3$: $e(G)=\frac{1}{2} pn^2 \pm o(n^2)$ and $G$ has $\frac{1}{8}p^4n^4 \pm o(n^4)$ cycles of length $4$.
\end{thm}

Throughout this paper, unless specified otherwise,
when considering a subset of vertices $U \subset V$ such that $|U|=\alpha n$ for some $\alpha$,
we tacitly assume that $|U|=\lfloor \alpha n \rfloor$ or $|U|=\lceil \alpha n \rceil$.
Since we mostly consider asymptotic values, this difference will not affect our calculation.

For a positive real $\delta$, we say that a graph $G$ is
\emph{$\delta$-close to satisfying $\mathcal{P}_1$} if $e(U) =
\frac{1}{2}p|U|^2 \pm \delta n^2$ for all $U \subset V(G)$, and
similarly define it for other properties. The formal definition of
equivalence of properties in Theorem \ref{thm_chunggraham} is as
following: for every $\varepsilon > 0$, there exists a $\delta$ such
that if a graph is $\delta$-close to satisfying one property, then
it is $\varepsilon$-close to satisfying another.

We call a graph \emph{$p$-quasi-random}, or \emph{quasi-random} if
the density $p$ is clear from the context, if it satisfies
$\mathcal{P}_1$, and consequently satisfies all of the equivalent
properties of Theorem \ref{thm_chunggraham}. We also say that a
graph property is \emph{quasi-random} if it is equivalent to
$\mathcal{P}_1$. Note that the random graph $G(n,p)$ with high
probability is $p$-quasi-random. However, it is not true that all
the properties of random graphs are quasi-random. For example, it is
easy to check that the property of having $\frac{1}{2}pn^2 + o(n^2)$
edges is not quasi-random (as an instance, there can be many
isolated vertices). For more details on quasi-random graphs we refer
the reader to the survey of Krivelevich and Sudakov
\cite{kri_sudakov}. Quasi-randomness was also studied in many other
settings besides graphs, such as set systems \cite{quasiset},
tournaments \cite{quasitournament} and hypergraphs
\cite{quasihypergraph}.

The main objective of our paper is to study the quasi-randomness of
graph properties given by certain graph cuts. These kind of
properties were first studied by Chung and Graham in \cite{quasiset,
cgqrgraph}. For a real $\alpha \in (0,1)$, the cut property
$\PCUT{\alpha}$ is the collection of graphs $G$ satisfying the
following: for any $U \subset V(G)$ of size $|U|=\alpha n$, we have
$e(U, V \backslash U)=p\alpha(1-\alpha)n^2+o(n^2)$.
%
As it turns out, for most values of $\alpha$, the cut property
$\PCUT{\alpha}$ is quasi-random. In \cite{quasiset, cgqrgraph}, the
authors proved the following beautiful theorem which characterizes
the quasi-random cut properties.

\begin{thm} \label{thm_intro1}
$\PCUT{\alpha}$ is quasi-random if and only if $\alpha \neq 1/2$.
\end{thm}

To see that $\PCUT{1/2}$ is not quasi-random, Chung and Graham \cite{quasiset, cgqrgraph} 
observed that the graph obtained
by taking a random graph $G(n/2, 2p)$ on $n/2$ vertices and an independent set on the remaining $n/2$
vertices, and then connecting these two graphs with a random bipartite graph with edge probability
$p$, satisfies $\PCUT{1/2}$ but is not quasi-random.

A $r$-cut is a partition of a vertex set $V$ into subsets $V_1, \cdots,
V_r$, and if for a vector $\VALPHA=(\alpha_1, \cdots, \alpha_r)$,
the size of the sets satisfies $|V_i| = \alpha_i |V|$ for all $i$,
then we call this an $\VALPHA$-cut. An $\VALPHA$-cut is called
\emph{balanced} if $\VALPHA=(1/r, \cdots, 1/r)$ for some $r$, and is
\emph{unbalanced} otherwise. For a $k$-uniform hypergraph $G$ and a
cut $V_1, \cdots, V_r$ of its vertex set, let $e(V_1, \cdots, V_r)$
be the number of hyperedges which have at most one vertex in each
part $V_i$ for all $i$.

A $k$-uniform hypergraph $G$ is \emph{(weak) $p$-quasi-random} if for every
subset of vertices $U \subset V(G)$, $e(U)=p\frac{|U|^k}{k!} \pm
o(n^k)$. Let $\PCUT{\VALPHA}$ be the following property: for every
${\VALPHA}$-cut $V_1, \ldots, V_r$, $e(V_1, \cdots, V_r)=(p+o(1))n^k
\sum_{S \subset [r], |S|=k} \prod_{i \in S} \alpha_i$. 
Note that previously we mentioned the example which illustrate 
the non-quasi-randomness of $\PCUT{1/2}$. 
As noticed by Shapira and Yuster \cite{ShaYus}, a similar construction as above shows
that $\PCUT{1/k, \cdots, 1/k}$ is not quasi-random.
In fact, they generalized Theorem~\ref{thm_intro1} by proving
the following theorem.

\begin{thm} \label{thm_intro2}
Let $k \ge 2$ be a positive integer. For $k$-uniform hypergraphs,
the cut property $\PCUT{\VALPHA}$ is quasi-random if and only if
${\VALPHA} \neq (1/r, \ldots, 1/r)$ for some $r \ge k$.
\end{thm}

\medskip

For a fixed graph $H$, let $\mathcal{P}_H$ be the following property
: for every subset $U \subset V$, the number of copies of $H$ in $U$
is $(p^{|E(H)|}+o(1)){|U| \choose |V(H)|}$. In \cite{simonqrgraph},
Simonovits and S\'{o}s proved that $\mathcal{P}_H$ is equivalent to
$\mathcal{P}_1$ and hence is quasi-random. For a fixed graph $H$, as
a common generalization of Chung and Graham's and Simonovits and
S\'{o}s' theorems, we can consider the number of copies $H$ having
one vertex in each part of a cut. Let us consider the cases when $H$
is a clique of size $k$.

\begin{dfn} \label{dfn_cliquecut}
Let $k,r$ be positive integers such that $r \ge k \ge 2$, and let
${\VALPHA}=(\alpha_1,\cdots,\alpha_r)$ be a vector of positive real
numbers satisfying $\sum_{i=1}^{r} \alpha_i=1$. We say that a graph
satisfies the $K_k$ cut property $\mathcal{C}_{k}({\VALPHA})$ if for
every ${\VALPHA}$-cut $(V_1, \cdots, V_r)$, the number of copies of
$K_k$ which have at most one vertex in each of the sets $V_i$ is
$(p^{\binom{k}{2}} \pm o(1))n^k \sum_{S \subset [r], |S|=k} \prod_{i
\in S} \alpha_i$.
\end{dfn}

Shapira and Yuster \cite{ShaYus} proved that for $k \ge 3$,
$\mathcal{C}_{k}({\VALPHA})$ is quasi-random if ${\VALPHA}$ is
unbalanced (note that $\mathcal{C}_2(\VALPHA)$ is quasi-random if
and only if $\VALPHA$ is unbalanced). This result is a corollary of
Theorem \ref{thm_intro2} by the following argument. For a graph $G$
satisfying $\mathcal{C}_{k}(\VALPHA)$, consider the $k$-uniform
hypergraph $G'$ on the same vertex set where a $k$-tuple of vertices
forms an hyperedge if and only if they form a clique in $G$. Then
$G'$ satisfies $\PCUT{\VALPHA}$ and thus is quasi-random. By the
definition of the quasi-randomness of hypergraphs, this in turn
implies that the number of cliques of size $k$ inside every subset
of $V(G)$ is ``correct'', and thus by Simonovits and S\'{o}s'
result, $G$ is quasi-random.

Note that for balanced $\VALPHA$ this approach does not give enough
information, since it is not clear if there exists a graph whose
hypergraph constructed by the above mentioned process is not
quasi-random but satisfies $\PCUT{\VALPHA}$ (nonetheless as the
reader might suspect, the properties $\PCUT{\VALPHA}$ and
$\mathcal{C}_{k}(\VALPHA)$ are closely related even for balanced
$\VALPHA$). Shapira and Yuster made this observation and left the
balanced case as an open question asking whether it is quasi-random
or not (in fact, they asked the question 
for $\VALPHA = (1/k, \cdots, 1/k)$, but here we consider the slightly 
more general question for all balanced $\VALPHA$ as mentioned above).
Janson \cite{Janson} independently posed the same question
in his paper that studied quasi-randomness under the framework of graph limits.
In this paper, we settle this
question by proving the following theorem :

\begin{thm} \label{main_intro}
Fix a real $p \in (0,1)$ and positive integers $r,k$ such that $r
\ge k \ge 3$. For every positive $\varepsilon$, there exists a positive
$\delta$ such that the following is true. If $G$ is a graph
which has density $p$ and is $\delta$-close to satisfying the $K_k$
balanced cut property $\mathcal{C}_{k}(1/r,\cdots,1/r)$, then $G$ is
$\varepsilon$-close to being $p$-quasi-random.
\end{thm}


The rest of the paper is organized as follows. In
Section \ref{notation} we introduce the notations
we are going to use throughout the paper and state previously
known results that we need later. In Section
\ref{main} we give a detailed proof of the most important
base case of Theorem \ref{main_intro}, triangle balanced cut property,
i.e. $\mathcal{C}_{3}(1/r,\cdots,1/r)$.
In Section \ref{general}, we prove the general case
as a consequence of the base case.
The last section contains some concluding remarks and open problems
for further study.

\section{Preliminaries} \label{notation}

Given a graph $G=(V,E)$ and two vertex sets $X,Y \subset V(G)$, we
denote by $E(X,Y)$ the set of edges which have one end point in $X$
and the other in $Y$. Also we write $e(X,Y)=|E(X,Y)|$ to indicate
the number of edges and $d(X,Y) = \frac{e(X,Y)}{|X||Y|}$ for the
density. For a cut ${\bf X} = (X_1,\cdots,X_r)$ of the vertex set,
we say that a triangle with vertices $u, v, w$ \emph{crosses the
cut} ${\bf X}$ if it contains at most one vertex from each set, and
denote it by $(u,v,w) \pitchfork {\bf X}$. We use $Tr( {\bf X} )$
for the number of triangles with vertices $(u,v,w) \pitchfork {\bf
X}$. For a $k$-uniform hypergraph and a partition $V_1, \ldots, V_t$
of its vertex set $V$ into $t$ parts, we define its \emph{density
vector} as the vector in $\mathbb{R}^{t \choose k}$ indexed by the
$k$-subsets of $[t]$ whose $\{i_1,\cdots,i_k\}$-entry is the density
of hyperedges which have exactly one vertex in each of the sets
$V_{i_1}, \cdots, V_{i_k}$. Throughout the paper, we always use
subscripts such as $\delta_{\ref{thm_ShaYus}}$ to indicate that
the parameter $\delta$ comes from Theorem \ref{thm_ShaYus}.

To state asymptotic results, we utilize the following standard
notations. For two positive-valued functions $f(n)$ and $g(n)$, write $f(n)=
\Omega(g(n))$ if there exists a positive constant $c$ such that
$\lim \inf_{n \rightarrow \infty} f(n)/g(n) \geq c$, $f(n) =
o(g(n))$ if $\lim \sup_{n \rightarrow \infty} f(n)/g(n) = 0$. Also,
$f(n) = O(g(n))$ if there exists a positive constant $C > 0$ such
that $\lim \sup_{n \rightarrow \infty} f(n)/g(n) \leq C$.

\medskip

To isolate the unnecessary complication arising from the error
terms, we will use the notation $x=_\varepsilon y$ if
$|x-y|=O(\varepsilon)$ and say that $x,y$ are $\varepsilon$-equal.
For two vectors, we define $\vec{x}=_\varepsilon \vec{y}$ if
$\|\vec{x}-\vec{y}\|_{\infty}=O(\varepsilon)$. We omit the proof of
the following simple properties (we implicitly assume that the
following operations are performed a constant number of times in
total). Let $C$ and $c$ be positive constants.

\noindent (1a)~(Finite transitivity)~If $x=_\varepsilon y$ and $y=_\varepsilon z$, then $x=_\varepsilon z$.\\
(1b)~(Complete transitivity)~For a finite set of numbers $\{x_1,
\cdots, x_n\}$. If $x_i=_\varepsilon x_j$ for every $i,j$, then
there exists $x$ such that $x_i=_\varepsilon x$ for all $i$.\\
(2)~(Additivity)~If $x=_\varepsilon z$ and $y=_\varepsilon w$, then $x+y=_\varepsilon z+w$.\\
(3)~(Scalar product)~If $x=_\varepsilon y$ and $0<c \leq a \leq C$, then $ax=_{\varepsilon} ay$ and $x/a=_\varepsilon y/a$.\\
(4)~(Product)~If $x,y,z,w$ are bounded above by $C$, then $x=_\varepsilon y$ and $z=_\varepsilon w$ implies that $xz=_\varepsilon yw$.\\
(5)~(Square root)~If both $x$ and $y$ are greater than $c$, then $x^2=_\varepsilon y^2$ implies that $x=_\varepsilon y$.\\
(6)~For the linear equation $A\vec{x}=_\varepsilon \vec{y}$, if all the entries of an invertible matrix $A$ are bounded by  $C$, and the determinant of $A$ is bounded from below by $c$, then $\vec{x}=_\varepsilon A^{-1}\vec{y}$.\\
(7)~If $xy=_\varepsilon 0$, then either $x=_{\sqrt{\varepsilon}} 0$ or $y=_{\sqrt{\varepsilon}} 0$.

\subsection{Extremal Graph Theory}

To prove the main theorem, we use the regularity lemma developed by
Szemer{\'e}di \cite{MR540024}. Let $G=(V,E)$ be a graph and
$\varepsilon > 0$ be fixed. A disjoint pair of sets $X,Y \subset V$
is called an \textit{$\varepsilon$-regular pair} if $\forall A
\subset X, B \subset Y$ such that $|A| \ge \varepsilon |X|, |B| \ge
\varepsilon |X|$ satisfies $|d(X,Y) - d(A,B)| \leq \varepsilon$. A
vertex partition $\{V_i\}_{i=1}^{t}$ is called an
\textit{$\varepsilon$-regular partition} if (i) the sizes of $V_i$
differ by at most 1, and (ii) $(V_i, V_j)$ is $\varepsilon$-regular for
all but at most $\varepsilon t^2$ pairs $1 \leq i < j \leq n$. The
regularity lemma states that every large enough graph admits a
regular partition. In our proof, we use a slightly different form
which can be found in \cite{MR1395865}:

\begin{thm}[Regularity Lemma] \label{thm_regularitylemma}
For every real $\varepsilon>0$ and positive integers $m,r$ there
exists constants $T(\varepsilon,m)$ and $N(\varepsilon,m)$ such that
given any $n \geq N(\varepsilon,m)$, the vertex set of any
$n$-vertex graph $G$ can be partitioned into $t$ sets $V_1, \cdots,
V_t$ for some $t$ divisible by $r$ and satisfying $m \leq t \leq
T(\varepsilon,m)$, so that
\begin{itemize}
\item $|V_i|<\lceil \varepsilon n\rceil$ for every $i$.
\item $||V_i|-|V_j|| \leq 1$ for all $i,j$.
\item Construct a reduced graph $H$ on $t$ vertices such that $i \sim j$ in $H$ if and only if
$(V_i,V_j)$ is $\varepsilon$-regular in $G$. Then the reduced graph has minimum degree at least $(1-\varepsilon)t$.
\end{itemize}
\end{thm}

As one can see in the following lemma, regular pairs are useful in
counting small subgraphs of a graph (this lemma can easily be
generalized to other subgraphs).

\begin{lemma} \label{lemma_counttriangles}
Let $V_1, V_2, V_3$ be subsets of vertices.
If the pair $(V_i, V_j)$ is $\varepsilon$-regular with density $d_{ij}$ for
every distinct $i,j$, then the number of triangles $Tr(V_1, V_2, V_3)$ is
\[ Tr(V_1, V_2, V_3) = (d_{12}d_{23}d_{31} + O(\varepsilon))|V_1||V_2||V_3|.\]
\end{lemma}
\begin{pf}
If a vertex $v \in V_1$ has degree $(1 + O(\varepsilon))d_{12}|V_2|$ in $V_2$ and
$(1 + O(\varepsilon))d_{13}|V_3|$ in $V_3$, then by the regularity of the pair
$(V_2, V_3)$, there will be $(1 + O(\varepsilon))|V_2||V_3|d_{12}d_{23}d_{31}$
triangles which contain the vertex $v$. By the regularity of the pair $(V_1, V_2)$, 
there are at least $(1-\varepsilon) |V_1|$ vertices in $V_1$ which have at least $(1 + O(\varepsilon))d_{12}|V_2|$ neighbors in $V_2$, 
and similar holds for the pair $(V_1, V_3)$.
Hence there are at least $(1 - 2\varepsilon)|V_1|$ such vertices satisfying both conditions.
Moreover, since each vertex in $V_1$ is contained in at most $|V_2||V_3|$ triangles,
there are at most $2\varepsilon |V_1||V_2||V_3|$ triangles which do
not contain such vertex from $V_1$. Therefore we have,
\[ Tr(V_1, V_2, V_3) = (1 + O(\varepsilon))|V_1||V_2||V_3|d_{12}d_{23}d_{31} + 2\varepsilon |V_1||V_2||V_3| = (d_{12}d_{23}d_{31} + O(\varepsilon))|V_1||V_2||V_3|.\]
\end{pf}

For a fixed graph $H$, a perfect $H$-factor of a large graph $G$ is
a collection of vertex disjoint copies of $H$ that cover all the
vertices of $G$. The next theorem is a classical theorem proved by
Hajnal and Szemer{\'e}di \cite{HajSze} which establishes a
sufficient minimum degree condition for the existence of a perfect
clique factor.

\begin{thm} [\cite{HajSze}] \label{thm_hajnalszemeredi}
Let $k$ be a fixed positive integer and $n$ be divisible by $k$. If $G$ is a graph on $n$ vertices
with minimum degree at least $(1 - 1/k)n$, then $G$ contains a perfect $K_k$-factor.
\end{thm}


\subsection{Concentration}

The following concentration result of Hoeffding \cite{MR0144363}
and Azuma \cite{azuma} will be used several times during
the proof (see also \cite[Theorem 3.10]{mcdiarmid}).

\begin{thm} [Hoeffding-Azuma Inequality] \label{thm_azuma}
Let $c_1, \ldots, c_n$ be constants, and let $X_1, \ldots, X_n$ be a
martingale difference sequence with $|X_k| \le c_k$ for each $k$. Then for any $t \ge 0$,
$$Pr \left( \left| \sum_{i=1}^n X_i \right| \geq t\right) \leq 2\exp\left(-\dfrac{t^2}{2\sum_{i=1}^n c_i^2}\right).$$
\end{thm}

The next lemma is a corollary of Hoeffding-Azuma's inequality.

\begin{lemma} \label{lemma_probabilisticlemma}
Let $G=(V,E)$ be a graph with $|V| = n$ and $|E| = d {n \choose 2}$
for some fixed real $d$. Let $U$ be a random subset of $V$ constructed by
selecting every vertex independently with probability $\alpha$. Then
$e(U) = \alpha^2 d {n \choose 2} + o(n^2)$ with probability at least
$1 - e^{-O(n^{1/2})}$.
\end{lemma}
\begin{pf}
Arbitrarily label the vertices by $1,\ldots, n$ and consider the
vertex exposure martingale. More precisely, let $X_k$ be the number of edges within $U$
incident to $k$ among the vertices $1, \ldots, k-1$ ($X_k = 0$ if $k \notin U$),
and note that $e(U) = X_1+\cdots+X_n$. Also note
that $\big( X_1+\cdots+X_k - \mathbb{E}[X_1+\cdots+X_k] \big)_{k=1}^n$ forms a martingale
such that $| X_k - \mathbb{E}[X_k] | \le n$ for all $k$.
Thus by Hoeffding-Azuma's inequality (Theorem \ref{thm_azuma}),
\[ Pr(|e(U) - \mathbb{E}[e(U)]| \ge C) \le 2e^{-2C^2/n^3}.  \]
Since $\mathbb{E}[e(U)] = \alpha^2 d {n \choose 2}$, by selecting $C = n^{7/4}$,
we obtain $e(U) = \alpha^2 d {n \choose 2} + o(n^2)$ with probability at least $1 - e^{-O(n^{1/2})}$
(see, e.g., \cite[Theorem 7.2.3]{AlSp} for more on vertex exposure martingales).
\end{pf}

Note that the probability of success in this lemma can be improved
by carefully choosing our parameters. However, Lemma \ref{lemma_probabilisticlemma} 
as stated is already strong enough for our later applications.

\subsection{Quasi-randomness of hypergraph cut properties}

Recall the cut property $\PCUT{\VALPHA}$ defined in the
introduction, and the fact that it is closely related to the clique
cut property $\mathcal{C}_k(\VALPHA)$. While proving Theorem
\ref{thm_intro2}, Shapira and Yuster also characterized the
structure of hypergraphs which satisfy the balanced cut property
$\PCUT{1/r,\cdots,1/r}$. Let $p \in (0,1)$ be fixed and $t$ be an
integer. In order to classify the $k$-uniform hypergraphs satisfying
the balanced cut property, we first look at certain edge-weighted
hypergraphs.
Fix a set $I \subset [t]$ of size $|I| = t/2$, and consider the
weighted hypergraph on the vertex set $[t]$ such that the hyperedge
$e$ has density $2p|e \cap I|/k$ for all $e$.
Let ${\bf u}_{t,p,I}$ be the vector in $\mathbb{R}^{t \choose k}$ 
representing this weighted
hypergraph (each coordinate corresponds to a $k$-subset of $[t]$, and
the value of the vector at the coordinate is the edge-weight of that
hyperedge), and let $W_{t,p}$ be the affine subspace of
$\mathbb{R}^{t \choose k}$ spanned by the vectors
${\bf u}_{t,p,I}$ for all possible sets $I$ of size
$|I| = t/2$. In \cite{ShaYus}, the authors proved that the structure
of a (non-weighted) hypergraph which is $\delta$-close to satisfying
the balanced cut property $\PCUT{1/r,\cdots,1/r}$ can be described
by the vector space $W_{t,p}$ (note that the vector which has
constant weight lies in this space).

\begin{thm} [\cite{ShaYus}] \label{thm_ShaYus}
Let $p \in (0,1)$ be fixed. There exists a real $t_0$ such that for
every $\varepsilon > 0$, and for every $t \ge t_0$ divisible by $2r$
\footnote{The authors omitted the divisibility condition in their paper \cite{ShaYus}.},
there exists
$\delta = \delta(t, \varepsilon) > 0$ so that the following holds.
If $G$ is a $k$-uniform hypergraph with density $p$ which is
$\delta$-close to satisfying the balanced cut property
$\PCUT{1/r,\cdots,1/r}$, then for any partition of $V(G)$ into $t$
equal parts, the density vector ${\bf d}$ of this partition
satisfies $\|{\bf d}-{\bf y}\|_{\infty} \le \varepsilon$ for some
vector ${\bf y} \in W_{t,p}$.
\end{thm}

A part of the proof of Shapira and Yuster's theorem relies on
showing that certain matrices have full rank, and they establish
this result by using the following famous result from algebraic
combinatorics proved by Gottlieb \cite{gottlieb}. For a finite set
$T$ and integers $h$ and $k$ satisfying $|T| > h \ge k \ge 2$,
denote by $B(T,h,k)$ the $h$ versus $k$ inclusion matrix of $T$
which is the ${|T| \choose h} \times {|T| \choose k}$ $0$-$1$ matrix
whose rows are indexed by the $h$-element subsets of $T$, columns
are indexed by the $k$-elements subsets of $T$, and entry $(I,J)$ is
1 if and only if $J \subset I$.
\begin{thm} \label{thm_gottlieb}
$rank(B(T,h,k)) = {|T| \choose k}$ for all $|T| \ge h + k$.
\end{thm}

\section{Base case - Triangle Balanced Cut} \label{main}

In this section we prove a special case, triangle balanced cut
property, of the main theorem. Our proof consists of several steps.
Let $G$ be a graph which satisfies the triangle balanced cut
property. First we apply the regularity lemma to describe the
structure of $G$ by an $\varepsilon$-regular partition
$\{V_i\}_{i=1}^{t}$. This step allows us to count the edges or
triangles effectively using regularity of the pairs. From this point
on, we focus only on the cuts whose parts consist of a union of the
sets $V_i$. In the next step, we swap some vertices of $V_i$ and
$V_j$. By the triangle cut property, we can obtain an algebraic
relation of the densities inside $V_i$ and between $V_i$ and $V_j$.
After doing this, the problem is transformed into solving a system
of nonlinear equations, which basically implies that inside any
clique of the reduced graph, most of the densities are very close to
each other. Finally resorting to results from extremal graph theory,
we can conclude that almost all the densities are equal and thus
prove the quasi-randomness of triangle balanced cut property.

To show that our given graph is quasi-randomn, ideally, we would like to
show that the densities of edges between pair of parts in the regular 
partition is (almost) equal to each other. However, instead of
establishing quasi-randomness through verifying this strong condition, we will
derive it from a slightly weaker condition. 
More specifically, we will use the fact that if in an $\varepsilon$-regular
partition of the graph, the density of edges in most of the pairs of parts are
equal to each other, then the graph 
is quasi-random (there are some dependencies in parameters). Following is
the main theorem of this section.


\begin{thm}\label{maintheorem}
Fix a real $p \in (0,1)$ and an integer $r \ge 3$. For every
positive $\varepsilon$, there exists a positive real $\delta$ such that the
following is true. If $G$ is a graph which has density $p$ and is
$\delta$-close to satisfying the triangle balanced cut property
$\mathcal{C}_3(1/r,\cdots,1/r)$, then $G$ is $\varepsilon$-close to
being $p$-quasi-random.
\end{thm}
Let $G$ be a graph $\delta$-close to satisfying
$\mathcal{C}_3(1/r,\cdots,1/r)$. By applying the regularity lemma,
Theorem \ref{thm_regularitylemma}, to $G$, we get an
$\varepsilon$-regular equipartition $\pi = \{V_i\}_{i=1}^{t}$. We can
assume that $|V_1|=\cdots=|V_t|$ by deleting at most $t$ vertices.
The reason this can be done is that later when we use the triangle
cut property to count the number of triangles, the error term that
this deletion creates is at most $tn^2$ which is negligible
comparing to $\delta n^3$ when $n$ is sufficiently large. Also in
the definition of quasi-randomness, the error term from counting
edges is at most $tn$, which is also $o(n^2)$.

Now denote the edge density within $V_i$ by $x_i$, the edge density
between $V_i$ and $V_j$ by $d_{ij}$, and the density of triangles in
the tripartite graph formed by $(V_i,V_j,V_k)$ by $d_{ijk}$. Call a
triple $(V_i,V_j,V_k)$ {\em regular} if each of the three pairs is
regular.


Consider a family $\{\pi_{\alpha}\}_{\alpha \in [0,1]}$ of
partitions of $G$ given as follows:
$$\pi_{\alpha}=((1-\alpha)V_1+\alpha V_2, \alpha V_1+(1-\alpha)V_2,V_3,\cdots V_t).$$
In other words, we pick $U_1$ and $U_2$ both containing
$\alpha$-proportion of vertices in $V_1$ and $V_2$ uniformly at
random and exchange them to form a new equipartition $\pi_{\alpha}$.
To be precise, for fixed $\alpha$, the notation $\pi_\alpha$ 
represents a family of random partitions and not necessarily an individual partition.
For convenience we assume that
$\pi_\alpha$ is a partition constructed as above which satisfies
some explicit properties that we soon mention which a.a.s.~hold for
random partitions. Denote the new triangle density vector of
$\pi_{\alpha}$ by ${\bf d}^{\alpha}=(d_{ijk}^{\alpha})$.

Note that every $(1/r,\cdots,1/r)$-cut ${\bf X} = (X_1, \cdots,
X_r)$ of the index set $[t]$ also gives a $(1/r,\cdots,1/r)$-cut of
$V(G)$. With a slight abuse of notation, we use $(i,j,k) \pitchfork {\bf X}$ to 
indicate that $V_i$, $V_j$ and $V_k$ completely belongs to different parts of the cut induced by ${\bf X}$.

By the triangle balanced cut property, for every $\alpha \in
[0,1]$,
\[ (p^3 \pm \delta)\left(\dfrac{n}{r}\right)^3\cdot {r \choose 3} = \sum_{(i,j,k) \pitchfork
{\bf X}} Tr(V_i, V_j, V_k) = \sum_{(i,j,k) \pitchfork {\bf X}}
d_{ijk}^{\alpha}\left(\dfrac{n}{t}\right)^3. \]
So $\sum_{(i,j,k)
\pitchfork {\bf X}} d_{ijk}^{\alpha}=(p^3 \pm
\delta){r \choose 3}\left(\frac{t}{r}\right)^3$. Let $M$ be the
$\binom{t}{t/r,\cdots,t/r} \times \binom{t}{3}$ $0-1$ matrix whose
rows are indexed by the $(1/r,\cdots,1/r)$-cuts of the vertex set
$[t]$ and columns are indexed by the triples $\binom{[t]}{3}$, where
the $({\bf X},(i,j,k))$-entry of $M$ is $1$ if and only if $(i,j,k)
\pitchfork {\bf X}$. The
observation above implies $M {\bf d}^{\alpha}=(p^3 \pm
\delta){r \choose 3}\left(\frac{t}{r}\right)^3 \cdot \bf{1}$ where $\bf{1}$ is
the all-one vector. Thus if we let ${\bf d}'={\bf
d}^{1/2}-\frac{1}{2}{\bf d}^0-\frac{1}{2} {\bf d}^1$, then $M {\bf
d}'=_{\delta t^3} \bf{0}$. From this equation we hope to get useful
information about the densities $x_i$ and $d_{ij}$. With the help of
the following lemma, we can compute the new densities
$d_{ijk}^{\alpha}$, and thus the modified density vector ${\bf
d}'$, in terms of the densities $x_i$ and $d_{ij}$.

\begin{lemma}\label{computevector}
Let $\varepsilon$ satisfy $0<\varepsilon < d_{ij}/2$ for every $i,j$ and assume that the graph $G$ is
large enough. Then for all $\alpha \in (\varepsilon, 1 - \varepsilon)$, there exists a choice of sets $U_1, U_2$ such that the following holds.\\
\\
\noindent $(1)$ $$d_{ijk}^{\alpha}=
\begin{cases}
d_{ijk} & \text{if $\{i,j,k\} \cap \{1,2\} = \emptyset$}\\
(1-\alpha)d_{1jk}+\alpha d_{2jk}+o(1) & \text{if $i=1$ and $2 \not \in \{j,k\}$}\\
\alpha d_{1jk} + (1-\alpha) d_{2jk}+o(1) & \text{if $i=2$ and $1 \not \in \{j,k\}$}\\
\textrm{see (2)} & \text{if $i=1$ and $j=2$}
\end{cases}.
$$
$(2)$ If $(V_1,V_2,V_k)$ is a regular triple, then
$$d_{12k}^{\alpha}=((1-\alpha)^2+\alpha^2)d_{12}d_{1k}d_{2k}+\alpha(1-\alpha)(x_1d_{1k}^2+x_2d_{2k}^2) + O(\varepsilon).$$

\noindent $(3)$ Let $d_{ijk}^{'}=d_{ijk}^{\alpha}-(1-\alpha)d_{ijk}^{0}-\alpha d_{ijk}^{1}$. Then
$$d_{ijk}^{'}=
\begin{cases}
0 & \text{if $\{i,j,k\} \cap \{1,2\} = \emptyset$}\\
o(1) & \text{if $i=1$ and $2 \not \in \{j,k\}$}\\
o(1) & \text{if $i=2$ and $1 \not \in \{j,k\}$}
\end{cases}.
$$
Moreover, for the case $i=1$ and $j=2$, if $(V_1,V_2,V_k)$ is a regular triple, then
$$d_{12k}' =
\alpha(1-\alpha)(x_1d_{1k}^2+x_2d_{2k}^2-2d_{12}d_{1k}d_{2k})+O(\varepsilon).
$$
\end{lemma}
\begin{pf}
Throughout the proof, we rely on the fact that 
some events hold with probability $1-o(1)$. Since there are
fixed number of events involved, without
further mentioning, we will assume that all the involved events
happen together at the same time.

\medskip 

\noindent $(1)$ The claim clearly holds for the cases $\alpha=0$ and $\alpha=1$.

For $\alpha \in (0,1)$, if $\{i,j,k\} \cap \{1,2\} = \emptyset$, then the density $d_{ijk}^{\alpha}$
is not affected by the swap of vertices in $V_1$ and $V_2$ so it remains the same with $d_{ijk}$. In the case that $\{i,j,k\} \cap \{1,2\} = \{1\}$,
without loss of generality we assume $i=1$ and $j,k \neq 2$. We also assume that there are $S_x$ triangles with a fixed vertex $x \in V_1 \cup V_2$
and two other vertices belonging to $V_j$ and $V_k$ respectively (note that $S_x \le |V_j||V_k|$).
After swapping subset $U_1 \subset V_1$ with $U_2 \subset V_2$ such that $|U_1|=|U_2|=\alpha|V_i|$,
we know that the number of triangles in triple $(((V_1 \cup U_2) \backslash U_1), V_j, V_k)$ changes by
$\sum_{u \in U_2} S_u - \sum_{u \in U_1} S_u$.

Assume $|V_i|=m$, instead of taking $\alpha m$ vertices uniformly at
random, take every vertex in $V_1$ (or $V_2$) independently with
probability $\alpha$. This gives random variables $X_i$ for $1 \le i
\le m$ having Bernoulli distribution with parameter $\alpha$. Let
$R=\sum_{i=1}^m X_i$ and $S=\sum_{i=1}^m X_i S_i$. These random
variables represent the number of vertices chosen for $U_1$, and the
number of triangles in the triple that contain these chosen
vertices, respectively. It is easy to see
$$Pr(R = \alpha m)= \alpha^{\alpha m}(1-\alpha)^{(1-\alpha)m}\binom{m}{\alpha m} \sim \Omega \left(\dfrac{1}{\sqrt{\alpha(1-\alpha)}}m^{-1/2} \right), $$
and by Hoeffding-Azuma's inequality (Theorem \ref{thm_azuma})
$$Pr(|S-\mathbb{E}S| \geq C) \leq 2 \exp\left(-\dfrac{C^2}{2\sum_{i=1}^m S_i^2}\right) \leq 2 \exp\left(-\dfrac{C^2}{2m^3}\right).$$
Let $C=m^2$, and the second probability decreases much faster than
the first probability, thus we know that conditioned on the event
$R=\alpha m$, $S$ is also concentrated at its expectation
$\mathbb{E}S=\sum_{i=1}^m \alpha S_x=\alpha d_{1jk}m^3$. From here
we know the number of triangles changes by
$$\sum_{u \in U_2} S_u - \sum_{u \in U_1} S_u = \alpha d_{2jk}m^3 - \alpha d_{1jk}m^3 + o(m^3).$$
Therefore the new density is
$$d^{\alpha}_{1jk}=d_{1jk}+(\sum_{u \in U_2} S_u - \sum_{u \in U_1} S_u)/m^3=(1-\alpha)d_{1jk}+\alpha d_{2jk}+o(1).$$
We can use a similar method to compute $d_{2jk}^{\alpha}$ when $1 \not \in \{j,k\}$.

\medskip

\noindent $(2)$ Let $U_1 \subset V_1$, $U_2 \subset V_2$ be as in (1), and let
$V_1' = (V_1 \setminus U_1) \cup U_2 $, $V_2' = (V_2 \setminus U_2)
\cup U_1$. Then we have the identity
\[ Tr(V_1', V_2', V_k) = Tr(U_1, U_2, V_k) + Tr(U_1, V_1 \setminus U_1, V_k)
+ Tr(V_2 \setminus U_2, U_2, V_k) + Tr(V_2 \setminus U_2, V_1 \setminus U_1, V_k). \]

Since $\alpha \in (\varepsilon, 1 - \varepsilon)$, the triples
$(U_1, U_2, V_k)$ and $(V_1 \setminus U_1, V_2 \setminus U_2, V_k)$
are regular. Thus by Lemma \ref{lemma_counttriangles},
\begin{align*}
Tr(U_1, U_2, V_k) &= (d_{12}d_{1k}d_{2k} + O(\varepsilon))|U_1||U_2||V_k| = (d_{12}d_{1k}d_{2k} + O(\varepsilon))\alpha^2 m^3 \qquad \textrm{and,} \\
Tr(V_1 \setminus U_1, V_2 \setminus U_2, V_k) &= (d_{12}d_{1k}d_{2k} + O(\varepsilon))(1 - \alpha)^2 m^3 .
\end{align*}

To compute $Tr(U_1, V_1 \setminus U_1, V_k)$, let $E_1^{k} \subset
E(V_1)$ be the collection of edges such that their end points have $(d_{1k} \pm
\varepsilon)^2 m$ common neighbors in $V_k$. By the regularity of
the pair $(V_1, V_k)$, there are at most $2\varepsilon m$ vertices
in $V_1$ which do not have $(d_{1k} \pm \varepsilon)m$ neighbors
in $V_k$, otherwise taking this set of vertices and $V_k$ will contradict the regularity. 
If $v$ is not such a vertex, then since $\varepsilon < d_{1k} - \varepsilon$ by 
the hypothesis $\varepsilon< d_{ij}/2$, 
by using the regularity, we see that there are at most $2\varepsilon m$
other vertices in $V_1$ which do not have $(d_{1k} \pm
\varepsilon)^2 m$ common neighbors with $v$. Consequently there are
at most $4 \varepsilon m^2$ edges inside $V_1$ which do not have
$(d_{1k} \pm \varepsilon)^2 m$ common neighbors inside $V_k$. We call these edges ``exceptional''. Recall that 
$x_1$ denotes the density of edges in $V_1$, thus
$|E_1^k| \ge x_1 {m \choose 2} - 4\varepsilon m^2$. By Lemma
\ref{lemma_probabilisticlemma} and the calculation from part (1)
there exists a choice of $U_1$ of size $\alpha m$ such that,
\begin{align*} |E_1^k(U_1, V_1 \setminus U_1)| &= |E_1^k| - |E_1^k(U_1)| - |E_1^k(V_1 \setminus U_1)| \\
 &= (1 - \alpha^2 - (1-\alpha)^2 + o(1))|E_1^k|  = \alpha (1 - \alpha)x_1 m^2 + O(\varepsilon)m^2. \end{align*}
for all $k$. Note that the number of triangles $Tr(U_1, V_1
\setminus U_1, V_k)$ can be computed by adding the number of
triangles containing the edges in $E_1^k$ and then the number of
triangles containing the ``exceptional'' edges (recall that
there are at most $O(\varepsilon) m^2$ of the such edges). 
The latter can be crudely bounded by
$O(\varepsilon) m^2 \cdot m \le O(\varepsilon) m^3$. Since each edge
in $E_1^k$ is contained in $(d_{1k} \pm \varepsilon)^2 m$ triangles
(within the triple $(V_1, V_2, V_k)$),
\[ Tr(U_1, V_1 \setminus U_1, V_k) = |E_1^k(U_1, V_1 \setminus U_1)| \cdot (d_{1k}^2 + O(\varepsilon)) m + O(\varepsilon) m^3 = \alpha (1-\alpha) x_1 d_{1k}^2 m^3 + O(\varepsilon)m^3. \]
Similarly we can show that there exists a choice of $U_2$ of size $\alpha m$ that gives
\[ Tr(U_2, V_2 \setminus U_2, V_k) = \alpha (1-\alpha) x_2 d_{2k}^2 m^3 + O(\varepsilon)m^3. \]
for all $k$. Combining all the results together, we can conclude the existence of
sets $U_1$, $U_2$ such that
\[ d_{12k}^{\alpha} = ((1 - \alpha)^2 + \alpha^2)d_{12}d_{1k}d_{2k} + \alpha (1-\alpha)(x_1 d_{1k}^2 + x_2 d_{2k}^2) + O(\varepsilon). \]

Part $(3)$ is just a straightforward computation from the definition of $d'_{ijk}$ and $(2)$. 
\end{pf}


\begin{lemma}\label{getequation}
Let $\varepsilon$ satisfy $0 < \varepsilon < d_{ij}/2$ for every $i,j$
and assume that the graph $G$ is large enough. If $(V_i,V_j,V_k)$ is
a regular triple, then $x_id_{ik}^2+x_jd_{jk}^2-2d_{ij}d_{ik}d_{jk}
=_{\delta t^3+\varepsilon} 0$.
\end{lemma}
\begin{pf}
As mentioned before Lemma \ref{computevector}, the vector ${\bf d}'={\bf
d}^{1/2}-\frac{1}{2}{\bf d}^0-\frac{1}{2} {\bf d}^1$ satisfies $M
{\bf d}'=_{\delta t^3} \bf{0}$. For an index $k \neq 1,2$, consider
a balanced partition ${\bf X}$ of the vertex set $[t]$ such that $1$
and $2$ lies in different parts, and let $Y$ be the union of the
parts which contains neither 1 nor 2. Then by Lemma
\ref{computevector} (3),
\[ 0 =_{\delta t^3} \sum_{(i,j,k) \pitchfork {\bf X}} d_{ijk}' = t^3 \cdot o(1) + \sum_{k \in Y} d_{12k}', \]
where $o(1)$ goes to $0$ as the number of vertices in the graph $G$
grows. Since $Y$ can be an arbitrary set of size $(r-2)t/r$ not
containing $1$ and $2$, this immediately implies that $d_{12k}'
=_{\delta t^3} 0$ for all $k$. Thus if $(V_1, V_2, V_k)$ is a
regular triple, then by Lemma \ref{computevector} (3),
$x_1d_{1k}^2+x_2d_{2k}^2-2d_{12}d_{1k}d_{2k}=_{\delta
t^3+\varepsilon} 0$. By symmetry, we can replace $1$ and $2$ by
arbitrary indices $i,j$.
\end{pf}
\medskip

Using Theorem \ref{thm_ShaYus} which characterizes the non quasi-random
hypergraphs satisfying the balanced cut property, we can prove the following lemma which
allows us to bound the densities from below.


\begin{lemma} \label{boundlemma}
There exists $t_0$ such that for fixed $p \in (0,1)$ and every $t
\ge t_0$ which is divisible by $2r$, there exist $c = c(p)$ and $\delta_0 = \delta_0(t, p) > 0$
so that the following holds for every $\delta \le \delta_0$. If $G$
is a graph with density $p$ which is $\delta$-close to satisfying
the triangle balanced cut property, then for any partition $\pi$ of
$V(G)$ into $t$ equal parts, the density vector ${\bf d} =
(d_{ij})_{i,j}$ satisfies $d_{ij} \ge c$ for all distinct $i,j \in
[t]$.
\end{lemma}
\begin{pf}
Let $t_0 = t_{\ref{thm_ShaYus}}$, $\varepsilon = p^3/8$, and for a given $t
\ge t_0$ divisible by $2r$, let $\delta_0 = \min\{\delta_{\ref{thm_ShaYus}}(t, \varepsilon),
p^3/10 \}$. Let $V = V(G)$, and let $G'$ be the hypergraph over the
vertex set $V$ such that $\{i,j,k\} \in E(G')$ if and only if
$i,j,k$ forms a triangle in the graph $G$. Let $\pi$ be an arbitrary
partition of $V$ into $t$ equal parts $V_1, \ldots, V_t$, and let
$(d_{ij})_{i,j}$ be the density vector of the graph $G$, and
$(d_{ijk})_{i,j,k}$ be the density vector of the hypergraph $G'$
with respect to $\pi$. It suffices to show the bound $d_{ij} \ge
p^3/10$ for every distinct $i,j \in [t]$. For simplicity we will
only verify it for $d_{12}$. Note that the number of triangles which
cross $V_1, V_2, V_k$ is at most $e(V_1,V_2) \cdot |V_k| =
(|V_1||V_2|d_{12}) \cdot |V_k|$, and thus $d_{12k} \le d_{12}$ for
all $k \ge 3$. Consequently, by summing it up over all choices of $k$, 
we obtain the following inequality
which will be crucial in our argument:
\begin{align} \label{eqn_boundinglemma}
\sum_{k=3}^{t} d_{12k} \le (t-2) \cdot d_{12}.
\end{align}

Since $G$ is $\delta$-close to satisfying the triangle balanced cut
property, we know that the density $q$ of triangles is at least $q
\ge p^3 - \delta$. By Theorem \ref{thm_ShaYus}, $(d_{ijk})_{i,j,k}$
is $\varepsilon$-equal to some vector in $W_{t,q}$. Recall that the vectors
in $W_{t,q}$ can be expressed as an affine combination of the
vectors ${\bf u}_{t,q,I}=(u^I_{ijk})_{i,j,k}$ for sets $I \subset
[t]$ of size $|I| = t/2$, and note that the following is true no
matter how we choose the set $I$ (recall that by definition
we have $u^{I}_{12k}=\frac{2q}{3}|I \cap \{1,2,k\}|$) :
\begin{align*}
\sum_{k=3}^{t} u^I_{12k} \ge \sum_{k \in I \backslash \{1,2\}} \frac{2q}{3} \ge
\left(\frac{t}{2}-2\right) \frac{2q}{3}.
\end{align*}
Since $(d_{ijk})_{i,j,k}$ is $\varepsilon$-equal to an affine
combination of these vectors and an arbitrary affine combination
still satisfies the inequality above, for large enough $t$ and 
sufficiently small $\varepsilon$ we have
\begin{align} \label{eqn_boundinglemma2}
\sum_{k=3}^{t} d_{12k} \ge \left(\frac{t}{2}-2\right)\frac{2q}{3} -
t \varepsilon \ge \frac{tq}{3} - \frac{4q}{3} - \frac{tq}{6} \ge
\frac{tq}{8}.
\end{align}

By combining (\ref{eqn_boundinglemma}) and
(\ref{eqn_boundinglemma2}), we obtain $d_{12} \ge q/8 \ge (p^3 -
\delta)/8 \ge p^3/10$. Similarly we can deduce $d_{ij} \ge p^3/10$
for all distinct $i,j \in [t]$.
\end{pf}

Since Lemma \ref{boundlemma} asserts that all the
pairwise densities $d_{ij}$ are bounded from below by some constant,
we are allowed to divide each side of an $\varepsilon$-equality
by $d_{ij}$. This turns out to be a crucial ingredient in
solving the equations given by Lemma \ref{getequation}.


\begin{lemma} \label{structraltheorem}
Given a positive real $c$ and an integer $n \ge 4$, if $x_i \ge 0$ for every $i$, 
$d_{ij} \ge c$ for every distinct $i,j \in [n]$, and 
$x_{i}d_{ik}^2+x_jd_{jk}^2-2d_{ij}d_{ik}d_{jk}=_{\varepsilon}0$ 
for every distinct $i,j,k \in [n]$, 
then there exists $s \in [n]$, $x,y>0$ such that for any distinct
$i,j \neq s$, we have $d_{ij}=_{\varepsilon} \sqrt{x}$ and, for any $i
\neq s$, $d_{is}=_\varepsilon \sqrt{y}$. Moreover, $x_i
=_{\varepsilon} \sqrt{x}$ if $i \neq s$ and $x_s=_\varepsilon
\frac{\sqrt{x}}{y}(2y-x)$ (see, figure \ref{figure1}).
\end{lemma}

\begin{pf}
Throughout the proof, we heavily rely on the properties of
$\varepsilon$-equality given in Section \ref{notation}.

First consider the case $n=4$. By taking
$(i,j,k)=(1,2,3),(2,3,1),(3,1,2)$ respectively, we get the following
system of equations:
\begin{equation} \label{eq_baseeq}
\begin{cases}
d_{13}^2 x_1+d_{23}^2 x_2=_\varepsilon 2d_{12}d_{13}d_{23}\\
d_{12}^2 x_1+d_{23}^2 x_3=_\varepsilon 2d_{12}d_{13}d_{23}\\
d_{12}^2 x_2+d_{13}^2 x_3=_\varepsilon 2d_{12}d_{13}d_{23}
\end{cases}.
\end{equation}
Considering this as a system of linear equations with unknowns
$x_1,x_2,x_3$, the determinant of the coefficient matrix becomes
$2d_{12}^2d_{13}^2d_{23}^2 \ge 2c^6$. Moreover, the coefficients in
the matrix are bounded from above by $1$. Therefore we can solve the
linear system by appealing to property $(6)$ of
$\varepsilon$-equality and get
\begin{equation} \label{eq_linsol}
\begin{cases}
x_1=_\varepsilon \dfrac{d_{23}}{d_{12}d_{13}}(d_{12}^2+d_{13}^2-d_{23}^2)\\
x_2=_\varepsilon \dfrac{d_{13}}{d_{12}d_{23}}(d_{12}^2+d_{23}^2-d_{13}^2)\\
x_3=_\varepsilon \dfrac{d_{12}}{d_{13}d_{23}}(d_{13}^2+d_{23}^2-d_{12}^2)
\end{cases}.
\end{equation}
Then
\begin{equation} \label{eq_flipornot}
\begin{split}
x_1x_2 &=_\varepsilon \dfrac{d_{23}}{d_{12}d_{13}}(d_{12}^2+d_{13}^2-d_{23}^2) \cdot \dfrac{d_{13}}{d_{12}d_{23}}(d_{12}^2+d_{23}^2-d_{13}^2)\\
       &=_\varepsilon \dfrac{1}{d_{12}^2}[d_{12}^4-(d_{13}^2-d_{23}^2)^2]\\
       &=_\varepsilon \dfrac{1}{d_{12}^2}[d_{12}^4-(d_{14}^2-d_{24}^2)^2].
\end{split}
\end{equation}
The last equation comes from repeating the same step for the system
of equations for indices $1,2$, and $4$. Equation
(\ref{eq_flipornot}) implies $d_{13}^2-d_{23}^2=_\varepsilon
\pm(d_{14}^2-d_{24}^2)$, and $d_{ik}^2-d_{jk}^2=_\varepsilon
\pm(d_{il}^2-d_{jl}^2)$ for all distinct $i,j,k,l$ in general.
Assume that there exists an assignment $\{i,j,k,l\}=\{1,2,3,4\}$
such that $d_{ik}^2-d_{jk}^2=_\varepsilon -(d_{il}^2-d_{jl}^2)
\neq_\varepsilon 0$, (we call such case as a ``flip''). Without loss
of generality let $d_{13}^2-d_{23}^2=_\varepsilon
-(d_{14}^2-d_{24}^2)$. By equation (\ref{eq_baseeq}), we know that
$x_1d_{12}^2+x_3d_{23}^2=_\varepsilon
2d_{12}d_{13}d_{23}=_\varepsilon x_2d_{12}^2+x_3d_{13}^2$, from
which we get
$$d_{12}^2(x_1-x_2)=_\varepsilon (d_{13}^2-d_{23}^2)x_3.$$
Replace the index $3$ by $4$ and we get
$$d_{12}^2(x_1-x_2)=_\varepsilon (d_{14}^2-d_{24}^2)x_4.$$
By the assumption on a ``flip'', by subtracting the two equalities
we get $x_3+x_4=_\varepsilon 0$, thus $x_3=_\varepsilon 0$ and
$x_4=_\varepsilon 0$ by their nonnegativity. This is impossible from the equation
$x_3d_{13}^2+x_4d_{14}^2=_\varepsilon 2d_{13}d_{14}d_{34}$ and the
fact $d_{ij} \ge c$. Therefore no flip exists and we have
\begin{equation} \label{eq_noflip}
d_{ik}^2-d_{il}^2=_\varepsilon d_{jk}^2-d_{jl}^2 \qquad \forall
\{i,j,k,l\} = \{1,2,3,4\}.
\end{equation}

Since $d_{ij}\ge c$, the sum of $d_{12}^2+d_{13}^2-d_{23}^2$,
$d_{12}^2+d_{23}^2-d_{13}^2$ and $d_{13}^2+d_{23}^2-d_{12}^2$ is
equal to $d_{12}^2+d_{13}^2+d_{23}^2 \geq 3c^2$. So at least one of
the terms is greater than $c^2$, without loss of generality we can
assume $d_{12}^2+d_{13}^2-d_{23}^2 \geq c^2$. Recall that
$x_1=_\varepsilon
\dfrac{d_{23}}{d_{13}d_{12}}(d_{12}^2+d_{13}^2-d_{23}^2)$. By
equation (\ref{eq_noflip}), we also have $x_1 =_\varepsilon
\dfrac{d_{24}}{d_{14}d_{12}}(d_{12}^2+d_{14}^2-d_{24}^2)=_\varepsilon
\dfrac{d_{24}}{d_{14}d_{12}}(d_{12}^2+d_{13}^2-d_{23}^2)$. Therefore
$$\dfrac{d_{23}}{d_{13}}=_{\varepsilon}\dfrac{d_{24}}{d_{14}}.$$
By appealing to the bound $d_{ij} \ge c$, we get $d_{23}d_{14}=_{\varepsilon}d_{24}d_{13}$\
and $d_{23}^2d_{14}^2=_{\varepsilon}d_{24}^2d_{13}^2$, which implies
\begin{eqnarray*}
(d_{13}^2-d_{23}^2)(d_{13}^2-d_{14}^2) = d_{13}^2(d_{13}^2-d_{23}^2)-d_{13}^2d_{14}^2+d_{23}^2d_{14}^2=_{\varepsilon} d_{13}^2(d_{14}^2-d_{24}^2)-d_{13}^2d_{14}^2+d_{13}^2d_{24}^2=0.
\end{eqnarray*}
So either $d_{13}^2=_{\sqrt{\varepsilon}}d_{14}^2$ or
$d_{13}^2=_{\sqrt{\varepsilon}}d_{23}^2$. Thus at this point we may
assume the existence of indices $i,j,k$ satisfying $d_{ik}^2
=_{\sqrt{\varepsilon}} d_{jk}^2$. Assume that $d_{13}^2
=_{\sqrt{\varepsilon}} d_{14}^2$ as the other case can be handled
identically.

So $d_{13}^2=_{\sqrt{\varepsilon}}x$,
$d_{14}^2=_{\sqrt{\varepsilon}} x$ for some $x$ and by equation
(\ref{eq_noflip}) we have $d_{23}^2=_{\sqrt{\varepsilon}} y$,
$d_{24}^2=_{\sqrt{\varepsilon}} y$ for some $y$. We let
$d_{34}^2=z$, and the equation $d_{14}^2-d_{34}^2=d_{12}^2-d_{32}^2$
given by (\ref{eq_noflip}) translates to $d_{12}^2=_\varepsilon
x+y-z$. Moreover, from equation (\ref{eq_linsol}) for indices
$\{1,3,4\}$ and $\{1,2,4\}$ we know that
$$x_1=_\varepsilon \dfrac{d_{34}}{d_{14}d_{13}} (d_{14}^2+d_{13}^2-d_{34}^2)=_\varepsilon \dfrac{d_{24}}{d_{14}d_{12}} (d_{14}^2+d_{12}^2-d_{24}^2).$$
If we plug all equalities for $d_{ij}$ into this identity, we get
$$(2x-z)\dfrac{\sqrt{z}}{x}=_{\sqrt{\varepsilon}} (2x-z)\dfrac{\sqrt{y}}{\sqrt{x(x+y-z)}}. $$
So either $\dfrac{\sqrt{z}}{x}=_{\varepsilon^{1/4}}
\dfrac{\sqrt{y}}{\sqrt{x(x+y-z)}}$ or $z=_{\varepsilon^{1/4}}2x$. In
the first case, by solving this equation we get either
$z=_{\varepsilon^{1/8}}x$ or $z=_{\varepsilon^{1/8}}y$ (before
multiplying each side of the equation by its denominators, one must
establish the fact that $x+y-z$ is bounded away from 0. This can be
done by first noting that equation $\sqrt{\frac{z}{xy}}
=_{\varepsilon^{1/4}} \sqrt{\frac{1}{x+y-z}}$ holds, and then realizing
that the left hand side is bounded from above). Both of the above
solutions gives us values for $x_i$'s and $d_{ij}$'s as claimed (see figure \ref{figure1}, for
the case $z=_{\varepsilon^{1/8}}x$).

\begin{figure}
\centering
\includegraphics{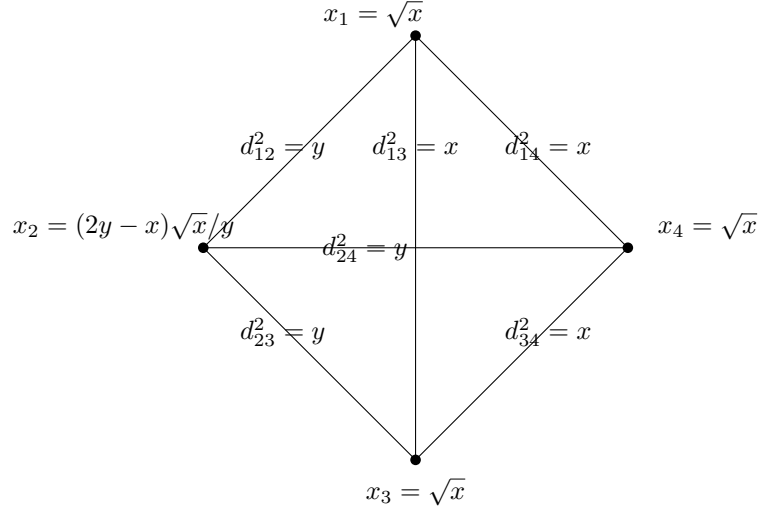}
\caption{Structure of solution for $n=4$ when $z=_{\varepsilon}x$, with vertices $1$ and $2$ permuted if $z=_\varepsilon y$.}
\label{figure1}
\end{figure}

In the second case $z=_{\varepsilon^{1/4}}2x$, we consider the
equation
$$x_2=_\varepsilon \dfrac{d_{34}}{d_{24}d_{23}} (d_{24}^2+d_{23}^2-d_{34}^2)=_\varepsilon \dfrac{d_{14}}{d_{24}d_{12}} (d_{24}^2+d_{12}^2-d_{14}^2)$$
to get
$$(2y-z)\dfrac{\sqrt{z}}{x}=_{\sqrt{\varepsilon}} (2y-z)\dfrac{\sqrt{y}}{\sqrt{x(x+y-z)}}. $$
By the previous analysis we may assume $z =_{\sqrt{\varepsilon}}
2y$, which implies $x=_{\varepsilon^{1/4}}y$ and
$d_{12}^2=_{\varepsilon^{1/4}} 0$. This is impossible by the fact
$d_{12} \ge c$.

Note that we have studied the case $n=4$. For $n=5$, suppose not all the edge densities are
$\varepsilon$-equal to the same value. In this case, there must be
four vertices such that not all the densities between them are equal.
Without loss of generality, $d_{12}=_\varepsilon
d_{13}=_\varepsilon d_{14}$ should be called $x$, and $d_{23}=_\varepsilon
d_{24}=_\varepsilon d_{34}$ should be called $y$, and $x \neq y$. Now let
us consider the collections of vertices $\{v_1,v_2,v_3,v_5\}$,
$\{v_1,v_2,v_4,v_5\}$, $\{v_1,v_3,v_4,v_5\}$. From the case $n=4$, we
know that $d_{15}=_\varepsilon x$, $d_{25}=_\varepsilon d_{35}
=_\varepsilon d_{45} =_\varepsilon y$. By repeating this process we can generalize it to arbitrary $n \geq 5$.
\end{pf}

Note that if in the regular partition, every pair of sets were
regular, then Lemma \ref{structraltheorem} itself forces the graph
to be quasi-random, as apart from one part (which is negligible),
all the densities are equal. However, the regularity lemma
inevitably produces a partition which contains some irregular pairs,
and in the remainder of the proof of Theorem \ref{maintheorem} we
will show how to handle this subtlety. The main idea is that since
there are only a small number of irregular pairs, the reduced graph
will contain many cliques, and thus that we can use Lemma
\ref{structraltheorem} to study its structure.

From now on in the reduced graph, when a clique of size at least $4$
is given, we will call the exceptional vertex $s$ ``bad'' and
all others ``good'' vertices. We also call the densities $x_s$
and $d_{is}$ for any $i \neq s$ ``bad'' and $d_{ij}$ ``good'' for
$i,j \neq s$. However as it will later turn out, most cliques of
size $4$ have $x=y$, and in this case we call every vertex and edge
``good''.

Now we can combine Lemmas \ref{getequation}, \ref{boundlemma}, and
\ref{structraltheorem} above to prove the main theorem which says
that the triangle balanced cut property is quasi-random.\\

\noindent \textbf{Proof of Theorem \ref{maintheorem} (triangle
case)}. Let $c = c_{\ref{boundlemma}}(p)$. We may assume that
$\varepsilon < \min \{ c/2, 1/32 \}$. Let $t_0 =
t_{\ref{boundlemma}}$ and
$T=T_{\ref{thm_regularitylemma}}(\varepsilon,t_0)$.  Let $\delta =
\min_{t_0 \le t \le T} \{ \varepsilon / t^3,
\delta_{\ref{boundlemma}}(t, p)\}$.

Let $G$ be a graph which is $\delta$-close to satisfying $\mathcal{C}_3(1/r,\cdots,1/r)$.
Consider the $\varepsilon$-regular equipartition $\pi$ of $G$:
$V(G)=V_1 \cup \cdots \cup V_t$ we mentioned before. This gives a
reduced graph $H$ on $t$ vertices of minimum degree at least
$(1-\varepsilon)t$ (we may assume that $t$ is divisible by $4r$). Every
edge $ij$ corresponds to an $\varepsilon$-regular pair $(V_i,V_j)$.
We mark on each edge of $H$ a weight $d_{ij}$ which is the density
of edges in $(V_i,V_j)$, and also the density $x_i$ inside $V_i$ on
the vertices. Parameters are chosen so that Lemma \ref{getequation}
and Lemma \ref{boundlemma} holds. Moreover, by the fact $\delta t^3
\le \varepsilon$, we have
$x_{i}d_{ik}^2+x_jd_{jk}^2-2d_{ij}d_{ik}d_{jk}=_{\varepsilon}0$ for
every regular triple $(V_i, V_j, V_k)$.
Thus whenever there is a clique of size at least $4$ in $H$, by
Lemma \ref{structraltheorem} we know that all the densities are
$\varepsilon$-equal to each other, except for at most one ``bad''
vertex. Since $\varepsilon < 1/32$ and $4 | t$, we can apply
Hajnal-Szemer\'{e}di theorem (Theorem \ref{thm_hajnalszemeredi}) to
the reduced graph $H$ and get an equitable partition of the vertices
of $H$ into vertex disjoint $4$-cliques $C_1, \cdots, C_{t/4}$.

For every $4$-clique $C_i$, from Lemma \ref{structraltheorem} we
know that there is at most one ``bad'' vertex. For two $4$-cliques
$C_i$ and $C_j$, we can consider the bipartite graph
$\mathcal{B}(C_i,C_j)$ between them which is induced from $H$. If
$\mathcal{B}(C_i,C_j)=K_{4,4}$, then it contains a subgraph
isomorphic to $K_{2,2}$ where all the vertices are ``good'' (two
vertices are good in $C_i$ and other two in $C_j$). If we apply the
structural lemma, Lemma \ref{structraltheorem}, to this new
$4$-clique (together with two edges coming from the two known
cliques), we get that the ``good'' densities of $C_i$ and $C_j$ are
$\varepsilon$-equal to each other.

Now consider the reduced graph $H'$ whose vertices correspond to the
$4$-cliques $C_i$, and $C_i$ and $C_j$ are adjacent in $H'$ if and
only if there is a complete bipartite graph between them. It is easy
to see that the minimum degree $\delta(H') \geq (1-16\varepsilon)|H'|$, since for any clique $C_i$, there are at most $4\varepsilon|H|$ edges intersecting 
it which is not in $H$, therefore $\delta(H') \geq |H'|-4\varepsilon|H| = (1-16\varepsilon)|H'|$. Take any two
vertices $u',v' \in V(H')$, since $d(u')+d(v') \geq 2(1-16\varepsilon)|H'|) >|H'|$ for $\varepsilon<1/32$, they have a
common neighbor $w'$, and thus by the discussion above, the ``good''
density in $C_{u'}$ or $C_{v'}$ are $\varepsilon$-equal to the
``good'' density in $C_{w'}$. So all the ``good'' densities are
$\varepsilon$-equal to each other. Thus by the total transitivity of
$\varepsilon$-equality (see, Section \ref{notation}), all the
``good'' densities are $\varepsilon$-equal to $p'$ for some $p'$.

We would like to show that $d_{ij} =_{\varepsilon} p'$ for all but
at most $O(\varepsilon) t^2$ edges of the reduced graph $H$. We
already verified this for ``good'' edges $\{i,j \}$ belonging to the
cliques $C_1, \ldots, C_{t/4}$. If $C_i$ is adjacent to $C_j$ in
$H'$ then they actually form a clique of size $8$ in $H$, and by Lemma
\ref{structraltheorem} there is at most one ``bad'' vertex there. Hence 
there is at most one ``bad vertex'' in two adjacent cliques. 
Consequently the total number of cliques that contain at least one
``bad'' vertex cannot exceed the independence number of $H'$, which
is at most $|H'|-\delta(H') \leq 16\varepsilon t$. Thus among the
cliques $C_1, \ldots, C_{t/4}$ there are at most $16\varepsilon t$
cliques which contain at least one ``bad'' vertex. Moreover, the
density of an edge in $H$ which is part of a $K_{4,4}$ connecting
two ``good'' cliques $C_i$ and $C_j$ are $\varepsilon$-equal to $p'$
again by Lemma \ref{structraltheorem}. Among the remaining edges,
all but at most $16\varepsilon t^2$ are such edges connecting two ``good'' cliques with
a $K_{4,4}$ as
otherwise $e(H) < {t \choose 2} - \varepsilon t^2$ which is a contradiction.
Therefore all but at most
$O(\varepsilon)t^2$ edges of $H$ have density $\varepsilon$-equal to
$p'$. This in turn implies that the density of $G$ is equal to $p' +
O(\varepsilon)$. On the other hand we know that the density is $p$,
thus $p' =_\varepsilon p$.

Now by verifying that $G$ satisfies $\mathcal{P}_{2}(1/2)$ (see
Theorem \ref{thm_chunggraham}), we will show that $G$ is
quasi-random. For an arbitrary subset $U \subset V(G)$ of size
$n/2$, let us compute the number of edges in $e(U)$ and estimate its
difference with the number of edges of a subset of size $n/2$ in
$G(n,p)$. Here we estimate the number of edges in each part $U \cap V_i$ and between two different parts respectively:
\begin{equation*}
\begin{split}
~&~~\left| \, e(U)-\binom{n/2}{2}p \, \right|\\
\leq & \sum_{i=1}^t \left| \, e(U \cap V_i)-\binom{|U \cap V_i|}{2}p \, \right|+\sum_{i,j} \left| \, 
e(U \cap V_i, U \cap V_j)-|U \cap V_i||U \cap V_j|p \, \right| ,\\
\end{split}
\end{equation*}
which by the fact that $e(U  \cap V_i, U \cap V_j) = |U \cap V_i| |U \cap V_j| (p'+O(\varepsilon))$ 
for all but at most $O(\varepsilon)t^2$ pairs $i,j$, is at most
\begin{equation*}
\begin{split}
& \sum_{i=1}^t |V_i|^2 + \left(\sum_{i,j} |U \cap V_i||U \cap V_j|\right)  \left(|p'-p|+O(\varepsilon)\right) + O(\varepsilon n^2) \\
\leq & n^2/t+O(\varepsilon |U|^2/2) + O(\varepsilon n^2) =  O(\varepsilon n^2).
\end{split}
\end{equation*}
In the last equation, we took $t$ to be sufficiently large depending on $\varepsilon$ and $p$. Therefore by the quasi-randomness of $\mathcal{P}_2(1/2)$,
we can conclude that $G$ is a quasi-random graph.
\hfill \qed


\section{General Cliques}
\label{general}

Throughout the section, $k$ and $r$ are fixed integers satisfying $r
\ge k \ge 4$. Let an $r$-balanced cut be a $(1/r, \cdots, 1/r)$-cut.
In this section, we will prove the remaining cases of the main
theorem, quasi-randomness of general $k$-clique  $r$-balanced cut
properties.

\begin{thm} \label{generalmaintheorem}
Fix a real $p \in (0,1)$ and positive integers $r,k$ such that $r
\ge k \ge 4$. For every $\varepsilon>0$, there exists a positive
real $\delta$ such that the following is true. If $G$ is a graph
which has density $p$ and is $\delta$-close to satisfying the $K_k$
balanced cut property $\mathcal{C}_{k}(1/r,\cdots,1/r)$, then $G$ is
$\varepsilon$-close to being $p$-quasi-random.
\end{thm}

Let $G$ be a graph which is $\delta$-close to satisfying the
$k$-clique $r$-balanced cut property. Apply the regularity lemma
(Theorem \ref{thm_regularitylemma}) to this graph to obtain an
$\varepsilon$-regular partition $\{V_i\}_{i=1}^{t}$ of the vertex set.
For $i \in [t]$, let $x_i$ be the density of the edges within $V_i$,
and for distinct $i,j \in [t]$, let $d_{ij}$ be the density of the
pair $(V_i, V_j)$. For $k \ge 2$, a $k$-tuple $J=\{i_1, \ldots, i_k\}$
is a multiset of $k$-indices (not necessarily distinct). Let $d_J$
be the density of $k$-cliques which have exactly one vertex in
each of the $V_{i_a}$ for $a=1,\ldots,k$. A $k$-tuple $J$ is called
\emph{regular} if $(V_{i_a}, V_{i_b})$ forms an
$\varepsilon$-regular pair for all $a,b \in [k]$. For a $k$-tuple
$J$ and a cut ${\bf X} = \{X_1, \ldots, X_r\}$, we say that $J$
crosses the cut ${\bf X}$ if $|J \cap X_i| \le 1$ for all $i$, and
denote it by $J \pitchfork {\bf X}$.

The proof of the $k$-clique $r$-balanced cut case follows the same
line of the proof of the triangle case. First we develop two lemmas
which correspond to Lemma \ref{getequation} and Lemma
\ref{boundlemma}. The following lemma can be regarded as generalization of 
Lemma \ref{getequation} for arbitrary integers $k \ge 3$.

\begin{lemma} \label{lemma_general_equation}
For any $k \geq 3$, let $\varepsilon$ be small enough depending on the densities $d_{ij}$ for all $i,j \in [t]$.
There exists a function $f: \mathbb{R} \rightarrow
\mathbb{R}$ such that the following holds. Let $J$ be a regular
$k$-tuple, $J' \subset J$ be such that $|J'| = k-2$, and $\{j_1,
j_2\} = J \setminus J'$. Then,
\[ x_{j_1} \left(\prod_{a \in J'} d_{a j_1 }\right)^2 \left(\prod_{a,b \in J', a < b} d_{a b}\right) + x_{j_2} \left(\prod_{a \in J'} d_{a j_2 }\right)^2 \left(\prod_{a,b \in J', a < b} d_{a b}\right) - 2\left(\prod_{a,b \in J, a < b} d_{a b}\right) =_{\varepsilon + \delta \cdot  f(t)} 0. \]
\end{lemma}
\begin{pf}
For the sake of clarity, without loss of generality we consider the case $\{1,2\} \subset J$ and $j_1 = 1, j_2 = 2$.
As in the triangle case, by considering the family of $t$-partitions
\[ \pi_{\alpha} = ((1-\alpha)V_1 + \alpha V_2, \alpha V_1 + (1 - \alpha)V_2, V_3, \ldots, V_t), \]
and the density vector $(d_{J}^{\alpha})_{J \in {[t] \choose k}}$
which arise from these partitions, we can define ${\bf d}' =
(d_{J}')_{J \in {[t] \choose k}}$ as $d_{J}' = d_{J}^{1/2} -
\frac{1}{2} d_J^0 - \frac{1}{2} d_J^1$. The same proof as in Lemma
\ref{computevector} gives us,
\begin{align}  \label{eqn_affinepart}
d_{J}' =
\left\{ \begin{array}{cll}
  & 0                   & \textrm{if $J \cap \{1,2\} = \emptyset$} \\
  & o(1)                                & \textrm{if $J \cap \{1,2\} = \{1 \}$} \\
  & o(1)                & \textrm{if $J \cap \{1,2\} = \{2 \}$}
  \end{array}\right.,
\end{align}
and if $\{1,2\} \subset J$ (let $J' = J \setminus \{1,2\})$) and $J$ is a regular
$k$-tuple, we have
\begin{align}  \label{eqn_nonaffine}
d_{J}' =  \alpha(1-\alpha) \left(\sum_{i=1}^2 x_{i} \left(\prod_{a \in J'} d_{a i }\right)^2 \left(\prod_{a,b \in J', a < b} d_{a b}\right)  - 2\left(\prod_{a,b \in J, a < b} d_{a b}\right)\right) + O(\varepsilon).
\end{align}
If $\{1,2\} \subset J$ and $J$ is not a regular $k$-tuple, then we do
not have any control on $d_{J}'$.

Let $M$ be the ${t \choose t/r,\ldots,t/r} \times {t \choose k}$
$0$-$1$ matrix whose rows are indexed by $r$-balanced cuts of the
vertex set $[t]$ and columns are indexed by the $k$-tuples ${[t]
\choose k}$. The $(\{X_1, X_2, \ldots, X_r\}, J)$-entry of $M$
is 1 if and only if $J \pitchfork \{X_1, X_2, \ldots, X_r\}$. We know
that $M {\bf d}' =_{\delta t^k} {\bf 0}$.

Consider the submatrix $N$ of $M$ formed by the rows of partitions
which have $1$ and $2$ in different parts, and columns of $k$-tuples
which include both $1$ and $2$. Let ${\bf d}''$ be the projection of
${\bf d}'$ onto the coordinates corresponding to the $k$-tuples
which contain both $1$ and $2$. By $M {\bf d}' =_{\delta t^k} {\bf 0}$, we can conclude that $N {\bf d}''
=_{\delta t^k}{\bf 0}$ given that the graph is large enough, since by \eqref{eqn_affinepart} only the columns $J$ with $\{1,2\} \subseteq J$ will play  a role 
here. Thus if we
can show that $N$ has full rank, then this implies that $d_J'
=_{\delta \cdot f(t)} 0$ for all $\{1,2\} \subset J$ (appeal to
property (6) of $\varepsilon$-equality - the entries of $N$ are
bounded and the size of $N$ depends on $t$).

Observe that the matrix $N$ has a lot of repeated rows 
(we obtain identical rows by swaping elements from 
the parts containing $1$ and $2$ with each other). However, 
these repetitions do not contribute to the rank and can be ignored.
Therefore we can consider a $0$-$1$ matrix $N'$
whose rows are indexed by the collection of subsets ${\bf Y} =
\{Y_1, \ldots, Y_{r-2}\}$ of the set $T= \{ 3, \ldots, t \}$ where
each part has size $t/r$, and columns are indexed by the
$(k-2)$-tuples $J' \in {T \choose k-2}$, where the entry $({\bf Y},
J')$ is 1 if and only if $J' \pitchfork {\bf Y}$. By fixing a subset
of $T$ of size $(r-2)t/r$ and considering all possible ${\bf Y}$
arising within this set, one can see that the row-space of $N'$ generates
the row-space of the $(r-2)t/r$ versus $k-2$ inclusion matrix of $T$, which we know
by Gottlieb's theorem, Theorem \ref{thm_gottlieb}, has full rank.
This implies that $N$ has full rank as well.
\end{pf}

Even though the equation which we obtained in Lemma
\ref{lemma_general_equation} looks a lot more complicated than the
triangle case, as it turns out, it is possible to make a
substitution of variables so that the equations above become exactly
the same as the equations in the triangle case. For a regular
$(k-3)$-tuple $I$, an index $j \notin I$, and distinct $j_1,j_2
\notin I$, define
\begin{align*} d_{j_1 j_2}^{I} &:= d_{j_1 j_2} \left( \prod_{a \in I} d_{a j_1} \right)^{1/2} \left(\prod_{a \in I} d_{a j_2} \right)^{1/2} \left( \prod_{a,b \in I, a < b} d_{a b} \right)^{1/3} \quad \textrm{and}  \\
x_{j}^{I} &:= x_j \left( \prod_{a \in I} d_{a j} \right) \left( \prod_{a,b \in I, a < b} d_{a b} \right)^{1/3}.
\end{align*}

\begin{claim} \label{claim_equationrelation}
Let $J$ be a regular $k$-tuple, $I \subset J$
be of size $|I| = k-3$, $\{ j_1, j_2, j_3\} = J \setminus I$, and $J' = I \cup \{j_3\}$. Then
\[ d_{j_1 j_2}^{I} d_{j_2 j_3}^{I} d_{j_3 j_1}^{I} = \prod_{a,b \in J, a < b} d_{a b}, \quad \textrm{and} \quad
 x_{j_1}^{I} (d_{j_1 j_3}^{I})^2 = x_{j_1} \left(\prod_{a \in J'} d_{a j_1 }\right)^2 \left(\prod_{a,b \in J', a < b} d_{a b}\right).\]
\end{claim}
\begin{pf} The claim follows from a direct calculation.
\end{pf}

In other words, Claim \ref{claim_equationrelation} transforms the
computation of the density of $K_r$ in the graph into the
computation of the density of triangles in another graph. This
observation will greatly simplify the equations obtained from Lemma
\ref{lemma_general_equation}.

\begin{lemma} \label{lemma_reducetotriangle}
Let $\varepsilon$ be small enough depending on the densities $d_{ij}$ for
all $i,j \in [t]$. There exists a function $f: \mathbb{R} \rightarrow
\mathbb{R}$ such that the following holds.  Let $J$ be a regular
$k$-tuple and $I \subset J$ be of size $|I| = k-3$. For $\{j_1, j_2,
j_3\} = J \setminus I$, we have $x_{j_1}^I (d_{j_1 j_3}^I)^2 +
x_{j_2}^I (d_{j_2 j_3}^I)^2 - 2 d_{j_1 j_2}^I d_{j_2 j_3}^I d_{j_3
j_1}^I =_{\varepsilon + \delta \cdot f(t)} 0$.
\end{lemma}
\begin{pf} This is an immediate corollary of Lemma \ref{lemma_general_equation} and Claim \ref{claim_equationrelation}.
\end{pf}

The next lemma corresponds to Lemma \ref{boundlemma} and establishes a lower bound on the densities.
We omit the proof which is a straightforward generalization of the proof of Lemma \ref{boundlemma}.

\begin{lemma} \label{lemma_general_lowerbound}
There exists $t_0$ such that for fixed $p \in (0,1)$ and every $t
\ge t_0$ divisible by $2r$, there exist $c = c(k,p)$ and $\delta_0 = \delta_0(t, p) >
0$ so that the following holds for every $\delta \le \delta_0$. If
$G$ is a graph with density $p$ which is $\delta$-close to
satisfying the $k$-clique balanced cut property, then for any
partition $\pi$ of $V(G)$ into $t$ equal parts, the density vector
${\bf d} = (d_{ij})_{i,j}$ satisfies $d_{ij} \ge c$ for all distinct
$i,j \in [t]$.
\end{lemma}

For every fixed regular $(k-3)$-tuple $I$, the set of equations that
Lemma \ref{lemma_reducetotriangle} gives is exactly the same as the
set of equations obtained from Lemma \ref{getequation}.
Consequently, by using Lemma \ref{lemma_general_lowerbound}, we can
solve these equations for every fixed $I$ just as in the triangle
case. Note that there is no need to (re)develop a statement corresponding
to Lemma \ref{structraltheorem}, since now that we reduced our problem
to the triangle case, the same lemma can be used as it is stated.

Thus as promised, we can reduce the case of general cliques to the
case of triangles. Therefore the proof of the triangle case of
Theorem \ref{maintheorem} can be repeated to give us useful information.
However, this observation does not immediately
imply that $d_{j_1 j_2} =_\varepsilon p$ for most of the pairs $j_1,
j_2 \not \in I$, since the only straightforward conclusion that we can draw is
that for every regular $(k-3)$-tuple $I$, there exists a constant
$p_{I}$ such that $d_{j_1 j_2}^{I} =_\varepsilon p_{I}$ for most of
the pairs $j_1, j_2 \notin I$. In order to prove the quasi-randomness of
balanced cut properties, we will need some control on the relation
between different $p_{I}$. Call a $k$-tuple $J$ \emph{excellent} if
it is regular, and for every $(k-3)$-tuple $I \subset J$, we have
$d_{j_1 j_2}^{I} =_\varepsilon p_I$ for all distinct $j_1, j_2 \in J
\setminus I$.

\begin{lemma} \label{lemma_excellenttuples}
Let $J$ be an excellent $k$-tuple. Then the density of every two pairs
in $J$ are $\varepsilon$-equal to each other.
\end{lemma}
\begin{pf}
For the sake of clarity, assume that $J = \{1,2, \ldots, k\}$. First,
consider $I = \{4,\ldots, k\}$. Then by the assumption, we have $d_{13}^I
=_\varepsilon d_{23}^I$, which by definition gives,
\[ d_{1 3} \left( \prod_{a \in I} d_{a 1} \right)^{1/2} \left(\prod_{a \in I} d_{a 3} \right)^{1/2} \left( \prod_{a,b \in I, a < b} d_{a b} \right)^{1/3} =_\varepsilon
d_{2 3} \left( \prod_{a \in I} d_{a 2} \right)^{1/2} \left(\prod_{a
\in I} d_{a 3} \right)^{1/2} \left( \prod_{a,b \in I, a < b} d_{a b}
\right)^{1/3}. \] After cancelation of the same terms, we can
rewrite this as,
\begin{align} \label{eqn_generalclique}
d_{1 3} \left(\prod_{a = 4}^{k} d_{a 1} \right)^{1/2}  =_\varepsilon
d_{2 3} \left(\prod_{a = 4}^{k} d_{a 2} \right)^{1/2} \quad  \Leftrightarrow  \quad
d_{1 3}^{1/2} \left(\prod_{a = 3}^{k} d_{a 1} \right)^{1/2}  =_\varepsilon
d_{2 3}^{1/2} \left(\prod_{a = 3}^{k} d_{a 2} \right)^{1/2}.
\end{align}
We can replace $3$ by $i$ for $i \in \{3, 4, \cdots, k\}$ and multiply each side of all these
equations to obtain,
\[ \prod_{i=3}^{k} \left( d_{1 i}^{1/2} \left(\prod_{a = 3}^{k} d_{a 1} \right)^{1/2}\right)  =_\varepsilon
\prod_{i=3}^{k} \left( d_{2 i}^{1/2} \left(\prod_{a = 3}^{k} d_{a 2} \right)^{1/2} \right), \]
which is equivalent to
\[ \left(\prod_{i = 3}^{k} d_{1 i} \right)^{(k-1)/2}  =_\varepsilon
\left(\prod_{i = 3}^{k} d_{2 i} \right)^{(k-1)/2}. \]
If we plug this back into equation (\ref{eqn_generalclique}), we get
$d_{13} =_\varepsilon d_{23}$. By repeating this process for other choice of indices,
we can conclude that the density of every two pairs are $\varepsilon$-equal to
each other.
\end{pf}

We now combine all these observations to show that $d_{e} =_\varepsilon p$ for
most of the edges $e$ of the reduced graph, which will in turn imply the quasi-randomness.

\medskip

\noindent \textbf{Proof of Theorem \ref{generalmaintheorem}}. Choose
$\varepsilon_0$ small enough depending on the constant $c =
c_{\ref{lemma_general_lowerbound}}(p)$ so that the condition of
Lemma \ref{lemma_reducetotriangle} holds, and let $f$ be the
function from Lemma \ref{lemma_reducetotriangle}. Let
$\varepsilon \le \min \{\varepsilon_0, 1/4\}$, $t_0 =
t_{\ref{lemma_general_lowerbound}}$, and let
$T=T_{\ref{thm_regularitylemma}}(\varepsilon,t_0)$. Let $\delta =
\min_{t_0 \le t \le T} \{\varepsilon / f(t),
\delta_{\ref{lemma_general_lowerbound}}(t,p)\}$.

Let $G$ be a graph which is $\delta$-close to satisfying the
$k$-clique $r$-balanced cut property. Apply the regularity lemma
(Theorem \ref{thm_regularitylemma}) to this graph to obtain an
$\varepsilon$-regular partition $\{V_i\}_{i=1}^{t}$ of the vertex set
where $t$ is divisible by $2r$.
For distinct $i,j \in [t]$, let $d_{ij}$ be the density of the pair
$(V_i, V_j)$. Note that the parameters are chosen so that Lemma
\ref{lemma_reducetotriangle} and Lemma
\ref{lemma_general_lowerbound} holds.

For every regular $(k-3)$-tuple $I$, define a graph $H_I$ as
following. The vertex set of $H_I$ is the collection of elements of
$[t] \setminus I$ which form a regular $(k-2)$-tuple together with
$I$. And $j_1, j_2 \in V(H_I)$ forms an edge if and only if the
$(k-1)$-tuple $I \cup \{j_1, j_2\}$ is regular. Since each part of the
regular partition forms a regular pair with at least $(1 -
\varepsilon)t$ of the other parts, we know that the graph $H_I$ has
at least $(1 - k \varepsilon)t$ vertices and minimum degree at least
$(1 - 2k \varepsilon)t$. Thus by Lemma
\ref{lemma_reducetotriangle}, Lemma \ref{lemma_general_lowerbound},
Lemma \ref{structraltheorem} and the proof of Theorem
\ref{maintheorem}, we know that there exists a $p_I$ such that at
least $(1 - O(\varepsilon))$-proportion of the edges of $H_I$ have
density $\varepsilon$-equal to $p_I$.

Select $k$ indices $j_1, \ldots, j_k$ out of $[t]$ independently and
uniformly at random. With probability at least $1 - O(\varepsilon)$,
the $k$-tuple is regular. Moreover, with probability at least $1-
O(\varepsilon)$, $d_{j_1 j_2}^{\{j_4, \ldots, j_k\}} =_\varepsilon
p_{\{j_4, \ldots, j_k\}}$ and the same is true for other choices of
indices as well. Therefore by the union bound, the $k$-tuple $\{j_1,
\ldots, j_k\}$ is excellent with probability at least $1 -
O(\varepsilon)$. Equivalently, the number of excellent $k$-tuples is
at least $(1 - O(\varepsilon)){t \choose k}$.

Call a pair of indices in $[t]$ \emph{excellent} if it is contained
in at least $\frac{2}{3}{t \choose k-2}$ excellent $k$-tuples.
Assume that there are $\eta t^2$ non-excellent edges. Then the
number of non-excellent $k$-tuples are at least
\[ \eta t^2 \times \frac{1}{3}{t \choose k-2} / {k \choose 2} = \Omega(\eta) {t \choose k}. \]
Therefore, $\eta = O(\varepsilon)$ and there are at most
$O(\varepsilon) t^2$ non-excellent edges. We claim that all the
excellent edges are $\varepsilon$-equal to each other. Take two
excellent edges $e, f$. Since each of these edges form an excellent
$k$-tuple with more than $\frac{2}{3}{t \choose k-2}$ of the
$(k-2)$-tuples, there exists a $(k-2)$-tuple which forms an
excellent $k$-tuple with both of these edges. Thus by Lemma
\ref{lemma_excellenttuples} applied to each of these $k$-tuples
separately, we can conclude that $d_{e} =_\varepsilon d_{f}$.

Consequently, by the total transitivity of $\varepsilon$-equality
(see, Section \ref{notation}), we can conclude that $d_{e} =_\varepsilon p'$
for some $p'$ for every excellent edge $e$. Then apply the same reasoning as in
the triangle case to show that $p' =_\varepsilon p$ and $G \in \mathcal{P}_2(1/2)$.
This proves the quasi-randomness of the graph $G$. \hfill \qed


\section{Concluding Remarks}

In this paper, we proved the quasi-randomness of $k$-clique balanced
cut properties for $k \geq 3$ and thus answered an open problem
raised by both Shapira-Yuster \cite{ShaYus} and Janson
\cite{Janson}. The most important base case was $k=3$ where we
solved a system of equations given by Lemma \ref{getequation}. The
existence of ``bad'' vertex in Lemma \ref{structraltheorem}
complicated the proof of the main theorem. It is hard to believe
that the case can be significantly simplified since even if we
assume that all the pairs are regular in the regular partition,
there is an assignment of variables $x_i$ and $d_{ij}$ which is not
all constant but forms a solution of the system.


\medskip

We conclude this paper with an open problem for further study. 
It is a generalization of balanced cut property to counting general graph $H$, which can also 
be regarded as an analogue of Simonovits and S\'{o}s' theorems for cuts. 

\begin{ques}
Let $k,r$ be positive integers satisfying $r \ge k \ge 3$. Let $H$
be a nonempty graph on $k$ vertices, and assume that every
$(1/r,\cdots,1/r)$-cut of a graph $G$ has the ``correct'' number of
copies of $H$ such that every vertex of $H$ is in a different part
of the cut. Does this condition force $G$ to be quasi-random?
\end{ques}

One might be able to adapt our approach to solve this question.
The main obstacle in this approach
lies in the fact that the new system of equations we get 
as in Lemma~\ref{getequation} now become
much more complicated to control.
In order to characterize the structure of 
densities, one will need to solve these system of equations and prove
statements such as in Lemma \ref{structraltheorem} and Lemma \ref{lemma_excellenttuples}.

\medskip


\noindent \textbf{Acknowledgement.} We would like to thank Asaf Shapira,
and our advisor Benny Sudakov for the kindness and advice they provided that greatly
helped us in doing this work. We also want to thank Svante Janson and the
two anonymous referees for their valuable comments and corrections.


\bibliographystyle{plain}

\end{document}